\documentclass[10pt]{amsart}
\usepackage{epsfig}
\usepackage[active]{srcltx}

\newtheorem{problem}{Problem}
\newtheorem{lemma}{Lemma}
\newtheorem{theorem}{Theorem}
\theoremstyle{remark}
\newtheorem{remark}{Remark}
\numberwithin{equation}{section}
\newcommand\be{\begin{eqnarray*}}
\newcommand\ee{\end{eqnarray*}}
\newcommand\ben{\begin{eqnarray}}
\newcommand\een{\end{eqnarray}}

\newcommand{\R}{\mathcal{R}}
\newcommand{\Rd}{\R^d}

\newcommand {\Yh}{{Y_h^*}}
\newcommand {\Ydiv}{{Y_{div}^*}}
\newcommand {\Q}{{Q}}
\newcommand {\Y}{{Y}}

\renewcommand {\u}{\mathbf{u}}
\renewcommand {\v}{\mathbf{v}}
\newcommand {\w}{\mathbf{w}}
\newcommand {\f}{\mathbf{f}}
\newcommand {\h}{\mathbf{h}}
\newcommand {\n}{\mathbf{n}}
\newcommand {\bp}{\mathbf{p}}
\newcommand {\bq}{\mathbf{q}}
\newcommand {\ttau}{\mathbf{\tau}}

\newcommand {\xxi}{\mathbf{\xi}}
\newcommand {\by}{\mathbf{Y}}

\newcommand{\eps}{\mathbb{\varepsilon}}
\def\sig{\mathbb{\sigma}}

\newcommand\C{\mathbb C}
\newcommand\bbL{{\mathbb L}\,}
\newcommand{\Mdsym}{\mathbb{ R}^{d \times d}_{sym}}
\newcommand{\cL}{\mathcal{L}}
\newcommand {\A}{\mathbb{A}}
\newcommand {\B}{\mathbb{B}}
\newcommand {\y}{\mathbb{Y}}
\newcommand {\M}{\mathbb{M}}

\def\IntO{\int\limits_\Omega}
\def\dx{\,dx}
\def\Norms#1#2{\left\Vert#1\right\Vert_{#2}}

\newcommand{\half}{\frac{1}{2}}
\newcommand {\Hdiv}{H(\mathrm{div};\Omega)}

\def\Norm#1{\left\Vert #1 \right\Vert}

\DeclareMathOperator{\diver}{div}

\def\pg{p_\Gamma}
\def\p{{\mathsf {p}}}
\def\q{{\mathsf {q}}}

\renewcommand {\ae}{\mbox{ a.e. in } \Omega}
\begin{document}

\title[A posteriori estimates for Barenblatt-Biot model]
{A posteriori error estimates for approximate solutions of
Barenblatt-Biot poroelastic model}

\author{J. M. Nordbotten}
\address{Department of Mathematics, University of Bergen, Norway}
\email{jan.nordbotten@math.uib.no}
\author{T. Rahman}
\address{Faculty of Engineering, Bergen University College, Norway}
\email{talal.rahman@hib.no}
\author{S. I. Repin }
\address{ V.A. Steklov Institute of Mathematics in St.-Petersburg,
Russia}
\email{repin@pdmi.ras.ru}
\author{J. Valdman}
\address{School of Engineering and Natural Sciences, University of Iceland, Iceland}
\email{janv@hi.is (corresponding author)}
\subjclass{Primary}

\begin{abstract}
The paper is concerned with the Barenblatt-Biott model in the theory of
poroelasticity. We derive a guaranteed estimate of the difference between exact
and approximate solutions expressed in a combined norm that encompasses
errors for the pressure fields computed from the diffusion part of the
model and errors related to stresses (strains) of the elastic part. Estimates do
not contain generic (mesh-dependent) constants and are valid
for any conforming approximation of pressure and stress fields.
\end{abstract}

\maketitle

\section{Introduction}
The standard mathematical model for diffusive flow in an elastic porous media is
the Biot's diffusion-deformation model of poroelasticity \cite{Biot} based on coupling between
the pore-fluid potential and the solid stress fields. The basic constitutive equations relate the
total stress to both the effective stress given by the strain of the
structure and to the potential arising from the pore-fluid. The model
consists of a momentum balance equation combined with Hooke's law
for elastic deformation, and a continuity equations combined with
Darcy's law. Originally, Biot's model was designed for homogeneous
porous media or single porosity media. The representation of
porosity and permeability in naturally occurring materials often
requires several distinct spatial scales. As for instance, in
reservoir model, the presence of heterogeneities like highly
permeable channels has a significant impact on the flow properties
of reservoir rock. Two of more scales of permeability are usually
observed, which is also referred to as dual permeability models.

Studies suggest that even for single-phase flow in relatively simple porous media, such as
sandstone, the fluid flows through a very small portion of the pore
space, while a greater part of it remains stagnant. A connected
system of highly permeable channels, characterized by relatively
simple pore space geometry, provides fluid flow through the
reservoir. The remainder of the reservoir, characterized by tortuous
pores and pore throats, is significantly less permeable. The highly
permeable channel component of a reservoir is relatively small, and
the remainder of the reservoir contains most of the fluid. This
contrast leads to the dual medium model of reservoir rock,
originally proposed by Barenblatt et al. \cite{BaZhKo} in the rigid
case. According to this model, the fluid flow in matrix blocks is
local, and only the local exchange of fluid between individual
blocks and the surrounding high permeable channels is supported.
This model contains a system of two diffusion equations, one for
each component, coupled by a distributed exchange term that, in its
simplest form, is proportional to the difference in potential between
fluids in the two components.

A combination of the Barenblatt's double-diffusion approach and
Biot's diffusion-deformation theory leads to what we call the
Barenblatt-Biot poroelastic model representing double diffusion in
elastic porous media. It takes the form
\ben
 -\nabla \cdot (\bbL \eps(\u)) + \alpha_1 \nabla p_1 + \alpha_2 \nabla p_2 &=& \f(x,t), \nonumber \\
c_1 \dot{p}_1 - \nabla \cdot (k_1 \nabla p_1) + \alpha_1 \nabla \cdot \dot{\u} + \kappa (p_1 - p_2) &=& h_1(x,t), \label{system} \\
c_2 \dot{p}_2 - \nabla \cdot (k_2 \nabla p_2) + \alpha_2 \nabla \cdot \dot{\u} + \kappa (p_2 - p_1) &=& h_2(x,t), \nonumber
\een
$\u$ is the displacement of the solid skeleton and $p_1$ and
$p_2$ are the fluid potentials in the respective components.
With the vector gradient operator $\nabla$, the linear Green strain tensor $\eps(\cdot)$ writes
\begin{equation} \label{eq:LinearizedElasticStrain}
    \eps (u) := \frac{1}{2} \left( \nabla \u + (\nabla \u)^T \right).
\end{equation}
The fourth-order elastic stiffness tensor $\bbL$ defines a stress tensor $\sig$ using the Hook's law
$$\sig:=\bbL \eps(\u).$$
In general, the permeabilities $k_1$ and $k_2$ may be heterogeneous
and anisotropic tensors, which may be functions of the
deformation. Herein, we will neglect this dependence and only consider
constant, scalar and homogeneous permeabilities.
Constants $\alpha_1$ and $\alpha_2$ measure changes of porosities due to
an applied volumetric strain. Mathematical analysis of this model based on the theory of
implicit evolution equations in Hilbert spaces is elaborated in \cite{SB}.

We note that multiple continua models are applicable to several other
porous media problems. We mention two cases in
particular. Firstly, contaminant transport
experiments clearly indicate that the particle dispersion is
non-Fickian, as reviewed in \cite{Berk06}. This makes both dual and
multiple continua models of interest, with dual media approaches already
common in applications. The use of more than two flowing continua was
argued by Gwo et al. \cite{Gwo96}, while a single flowing continua
coupled to multiple non-flowing continua (traps) is also reviewed in
\cite{Berk06}. The second application is heat transfer in fractured
rocks, in particular related to modelling of geothermal heat
extraction. Here, the slow interaction of diffusive heat transfer in
rock has to be modelled together with fast fluid flow in fractures. The
state of the art approach is to use multiple continua, frequently as
many as four or more \cite{Pruess}.

Our focus in this paper is to derive guaranteed and computable
bounds of approximation errors the static Barenblatt-Biot
system
\ben
\label{eq:3}
 -\nabla \cdot (\bbL \eps(\u)) + \alpha_1 \nabla p_1 + \alpha_2 \nabla p_2 &=& \f(x), \nonumber \\
 - \nabla \cdot (k_1 \nabla p_1) + \kappa (p_1 - p_2) &=& h_1(x), \\
 - \nabla \cdot (k_2 \nabla p_2) + \kappa (p_2 - p_1) &=& h_2(x), \nonumber
\een
which is considered in
bounded connected domain $\Omega \subset \Rd$
with Lipschitz continuous boundary $\Gamma$.
There are various  boundary conditions motivated by hydrological applications,
among which four boundary conditions, applicable to different parts of the
boundary $\Gamma = \bigcup \Gamma_i$ represent the most typical
cases. \\
{\em 1. Saturated land surface,
$\Gamma_1$, with infiltration and evaporation} is modelled as
\ben
\label{eq:3.1a}
&\sig(\u)\n=0\qquad&{\rm (
normal\, stress\,free\,condition )},\\
\label{eq:3.1b}
&\psi_{\Gamma_1} = \n\cdot(-k_1\nabla p_1 - k_2\nabla p_2)
&{\rm (normal\, fluid\, flux)},
\een
where $\n$ is the unit outward normal vector and $\psi_{\Gamma_1}$ is a given function.
We complete this boundary condition by
specifying that the normal component of the potential gradients at the
boundary are equal
\ben
\label{eq:3.1c}
\n\cdot\nabla (p_1-p_2)=0.
\een
In the case of constant $k$ (that we consider in this paper), the
condition \eqref{eq:3.1b} reads
\ben
\label{eq:3.1d}
&&\psi_{\Gamma_1} = -(k_1+ k_2)\n\cdot\nabla p_1=-(k_1+ k_2)\n\cdot\nabla p_2,
\een
which is in fact a version of the Darcy law at the boundary.
\\
{\em 2. Boundary to sea with a constant fluid potential (we call this boundary $\Gamma_2$).}
 Here we also impose normal stress as in \eqref{eq:3.1a}.
However, the boundary conditions for the potentials are of the Dirichlet
type, i.e.,
\ben
\label{eq:3.2}
&&p_1=p_2=p_{\Gamma_2}.
\een
{\em 3. Internal boundary with known head ($\Gamma_3$).} This may
represent either a fixed potential pumping well or the potential at
some measurement point. We model this as a no displacement boundary
 with Dirichlet conditions for the potentials as at
$\Gamma_3$, i.e.,
\ben
\label{eq:3.3a}
&& \u=0,\\
\label{eq:3.3b}
&& p_1=p_2=p_{\Gamma_3}.
\een
{\em 4. Impermeable bedrock, $\Gamma_4$.} Here, we impose no
displacement (as for $\Gamma_3$), and zero normal flux (as for
Equation \eqref{eq:3.1b} with $\psi_{\Gamma_4} = 0$).

For the unique solvability of the diffusion problem
one has to assume that
\be
{\rm meas}(\Gamma_2\cup\Gamma_3)\not=\emptyset.
\ee

\section{Variational formulation of the double diffusion system}
Since the displacement $\u$ is only involved in the first equation of
system \eqref{eq:3}, a double-diffusion problem
\ben
 - \nabla \cdot (k_1 \nabla p_1) + \kappa (p_1 - p_2) = h_1(x), \label{diffusion1}\\
 - \nabla \cdot (k_2 \nabla p_2) + \kappa (p_2 - p_1) = h_2(x) \label{diffusion2}
\een
is studied separately. It describes the steady flow of slightly compressible
fluid in a general heterogeneous medium consisting of two components.
Henceforth, we consider this problem
with the Dirichlet boundary conditions $p_1=p_2=\pg$ on $\Gamma$.
Let $\bar p$ be a function with square summable coefficients
that satisfies this boundary condition. It is convenient
to rewrite the problem in terms of new functions
$$\p_1:=p_1-\bar p, \quad \p_2:=p_2-\bar p.$$
Then, a weak formulation of ~\eqref{diffusion1}-\eqref{diffusion2} leads to

\begin{problem}
\label{problem_equality}
Assume that $(h_1, h_2) \in L^2(\Omega,\R^2)$. Find $
\bp = (\p_1, \p_2) \in H^1_0(\Omega,\R^2),$
satisfying the system of variational equalities
\begin{equation}
\begin{split}
 \IntO  k_1 \nabla \p_1 \cdot\nabla \q_1 + \IntO \kappa (\p_1 - \p_2) \q_1  \dx =
 \IntO (h_1(x) \q_1- k_1 \nabla \bar \p \cdot\nabla \q_1 )\dx  \\
 \IntO  k_2 \nabla \p_2\cdot \nabla \q_2 +
 \IntO \kappa (\p_2 - \p_1) \q_2  \dx =
 \IntO (h_2(x) \q_2-k_2 \nabla \bar \p \cdot\nabla \q_2) \dx
 \label{diffusions_weak}
 \end{split}
\end{equation}
for all testing functions $\bq = (\q_1, \q_2) \in H^1_0(\Omega,\R^2)$.
\end{problem}
This problem can be represented in a general form (which also encompasses other, more complicated models
of porous media).
For this purpose, we introduce the spaces
\ben
\Q := H^1_0(\Omega,\R^2) , \quad \Y:=L^2(\Omega,\R^{2 d}),
\een
and the corresponding dual spaces
\ben
\Q^*: = H^{-1}(\Omega,\R^2), \quad \Y^*:=L^2(\Omega,\R^{2 d}).
\een
Hereafter $L_2$ norms of all functions in $\Omega$ are denoted by $\Norm{\cdot}_\Omega$.
Duality pairings of $(\Q, \Q^*)$ and $(\Y, \Y^*)$ are denoted by
$\left< \cdot, \cdot \right>$ and $\left< \left< \cdot, \cdot \right> \right>$, respectively.
Also, we introduce a bounded linear operator $\Lambda \in \cL (\Q, \Y)$ and its adjoint
operator $\Lambda^* \in \cL(\Y^*, \Q^*) $ by the relations
\ben
 \Lambda \bq:= (\nabla \q_1, \nabla \q_2), \quad \Lambda^* \y^*=  ( -\diver y^*_1, -\diver y^*_2)^T.
\label{lambda_dual}
\een
The operators $\Lambda$ and $\Lambda^*$ satisfy the relation representing integration by parts
$$\left< \left< \y^*, \Lambda \bq \right> \right> =
\left< \Lambda^* \y^*, \bq \right> \quad  \mbox{for all } \y^* \in \Y^*, \bq \in \Q,$$
which can be written componentwise as
\ben
\IntO  \left( \by^*_1 \cdot \nabla \q_1 + \by^*_2 \cdot \nabla \q_2 \right)  \dx =
- \IntO  \left( \q_1 \diver \by^*_1  + \q_2 \diver \by^*_2  \right)  \dx,
\een
where $\bq =(\q_1, q_2)$ and $\y^* = (\by^*_1, \by^*_2)$. Now Problem \ref{problem_equality} can be represented in the form: Find $\bp \in \Q$ such that the equality
\ben
a(\bp, \bq) = l(\bq) \label{abstract_equality}
\een
holds for all $\bq \in \Q$. The bilinear form $a(\cdot, \cdot)$ and the linear form $l(\cdot)$ are defined as
\be
a( \bp,  \bq)&:=& \IntO \left(  \Lambda \bp:(\A \Lambda \bq) +  \bp \cdot \B\bq \right) \dx, \\
l(\bq)&:=& \IntO (\h \cdot \bq-\C\Lambda\bq) \dx,
\ee
$\A$, $\B$ and $\C$ are matrices formed by material dependent constants $k_1, k_2, k_3$,
\be
\A: = \begin{pmatrix} k_1 & 0 \\
    0 & k_2
    \end{pmatrix}, \quad
\B: = \begin{pmatrix}
    \kappa & -\kappa \\
    -\kappa & \kappa
    \end{pmatrix}, \quad
    \C: = \begin{pmatrix}
    k_1\nabla \bar p & 0 \\
    0 & k_2\nabla \bar p
    \end{pmatrix}
\ee
and $\h$ is the right hand side vector
\be
\h: = \begin{pmatrix}
  h_1 \\
  h_2
\end{pmatrix}.
\ee

\begin{remark}
We note that the symmetric matrix $\A$ is a positive definite matrix iff $k_1$ and $k_2$ are positive (since $\A \xxi \cdot \xxi \geq \min \{k_1, k_2 \} \Norm{\xxi}^2$
for all $\xxi \in \Rd$). However, $\B$ is symmetric but only positive semi-definite in case of the positive parameter $\kappa$, and its one-dimensional kernel is generated by
the vector $(1, 1)^T$.\\
\end{remark}
\begin{remark}
If $\bar p$ is sufficiently regular
(so that $\Lambda^* \C$ belongs to $\Y^*$), then
\be
l(\bq)&:=& \IntO (\h \cdot \bq-\Lambda^*\C \bq) \dx=
\IntO \widehat\h \cdot \bq \dx,
\ee
where
\be
\widehat\h: = \begin{pmatrix}
    h_1-\diver k_1\nabla\bar p \\
    h_2-\diver k_2\nabla\bar p
    \end{pmatrix}.
    \ee
\end{remark}
It is easy to verify that \eqref{abstract_equality} is the necessary condition for the minimizer of the following
convex variational problem.
\begin{problem}\label{problem_minimization}
Find $\bp \in \Q$ satisfying
\ben
F(\bp) + G(\Lambda \bp) =
\inf_{\bq \in \Q} \{ F(\bq) + G(\Lambda \bq) \},
\label{minimization_problem}
\een
where
\ben
&&F: \Q \rightarrow \R, \quad \quad F(\bq):=
\half \IntO   \bq \cdot \B\bq \dx - l(\bq) \label{F_form},
\een
and
\ben
&&G: \Y \rightarrow \R, \quad \quad G(\Lambda \bq):= \half \IntO  \Lambda \bq: (\A\Lambda \bq) \dx \label{G_form}.
\een
\end{problem}

\begin{theorem} [existence of unique solution]
Assume that $k_1, k_2 > 0$ and $\kappa \geq 0$. Then, there exists a unique solution
$\bp \in \Q$ of Problem \ref{problem_minimization}, which also represents the solution of Problem \ref{problem_equality}.
\begin{proof}
Existence of the unique minimizer follows from known results in the calculus of variations.
Indeed, under the give assumptions, the functional $F(\cdot) + G(\Lambda \cdot)$ is strictly convex and coercive
in the reflexive space $\Q$.
\end{proof}
\end{theorem}

\section{A posteriori error estimate of the double diffusion system}
In this section, we derive guaranteed and directly computable
bounds of the difference between exact and approximate solutions. Our analysis is based
upon a posteriori error estimation methods suggested in \cite{NeRe,Re3}.
Following the chapters 6 and 7 in \cite{NeRe}, first we need to find explicit forms of  dual functionals
\begin{equation}
\begin{split}
&F^*: \Q^* \rightarrow \R, \quad \quad F^*( \Lambda^* \y^*) := \sup_{\bq \in \Q} \{ \left< \Lambda^* \y^*, \bq \right>- F(\bq) \},  \\
&G^*: \Y^* \rightarrow \R, \quad \quad G^*(\y^*) := \sup_{\Lambda \bq \in \Y} \{\left< \left< \y^*,
\Lambda \bq \right> \right> - G(\Lambda \bq) \}.
\end{split}
\label{dual_functionals}
\end{equation}
and the corresponding compound functionals
\begin{equation}
\begin{split}
&D_F: \Q \times \Q^* \rightarrow \R, \quad \quad  D_F(\bq,\Lambda^* \y^*):=F(\bq) + F^*(\Lambda^* \y^*) - \left< \Lambda^* \y^*, \bq \right>,\\
&D_G: \Y \times \Y^* \rightarrow \R, \quad \quad D_G(\Lambda \bq, \y^*):=G(\Lambda \bq)  + G^*(\y^*) - \left< \left<\y^*, \Lambda \bq \right> \right>.
\end{split}
\label{compound_functionals}
\end{equation}
By the the sum of $D_F$ and $D_G$, we obtain the functional error majorant
\ben M(\bq,\y^*) := D_F(\bq,\Lambda^* \y^*) + D_G(\Lambda \bq, \y^*) \label{majorant}, \een
which provides a guaranteed upper bound of the error:
\ben
\half a(\bp-\bq, \bp - \bq) \leq M(\bq,\y^*) \quad \mbox{for all } \y^* \in \Y^*.
\label{estimate}
\een
The majorant is fully computable and depends only on the approximation $\bq \in \Q$
and arbitrary variable $\y^* \in \Y^* $.
\begin{lemma}[dual functionals]
For $k_1, k_2 > 0$ and $\kappa > 0$, it holds
\ben
G^*(\y^*) &=& \half \IntO \A^{-1} \y^*  : \y^* \dx,
\label{G_dual} \\
F^*(\Lambda^* \y^*) &=& \left\{ \begin{array}{ll}
\frac{1}{4 \kappa} \IntO (\Lambda^* \y^* + \h)^2 \dx   \quad \mbox{if } \Lambda^* y^*_1 + h_1  + \Lambda^* y^*_2 + h_2 = 0, \cr
+ \infty \quad \mbox{otherwise}.
\end{array}\right.
\label{F_dual}
\een
\begin{proof} The derivation of $G^*(\y^*)$ is straightforward, see \cite{NeRe}. The singularity of the matrix $\B$ makes the computation of $F^*(\Lambda^* \y^*)$ more technical.
\be
F^*(\Lambda^* \y^*) &=& \sup_{\bq \in \Q} \{ \left< \bq, \Lambda^* \y^* \right> - F(\bq) \}  \\
         &\geq& \sup_{\bq \in \Q: \q_1=\q_2} \{ \left< \bq, \Lambda^* \y^* \right> - F(\bq) \} \\
         &=&    \sup_{\q_1 \in H^1_0(\Omega)} \{ \left< \q_1, \Lambda^* \by^*_1 + \Lambda^* \by^*_2 \right> - F(\q_1,\q_1) \} \\
         &=&    \sup_{\q_1 \in H^1_0(\Omega)} \{ \left< \q_1, \Lambda^* \by^*_1 + h_1 + \Lambda^* \by^*_2  + h_2 \right> \} \\
         &=&    \left \{ \begin{array}{ll}
                    0 \quad \quad \mbox{if } \Lambda^* \by^*_1 + h_1 +  \Lambda^* \by^*_2  + h_2 = 0, \cr
                    + \infty \quad \mbox{otherwise}.
                         \end{array} \right.
\ee
Thus, finite values $F^*(\Lambda^* \y^*)$ are attained only on the subspace
\ben
\Lambda^* \by^*_1 + h_1 + \Lambda^* \by^*_2  + h_2 = 0 \label{subspace}
\een
and we must specially consider this case.
It holds
\be
F^*(\Lambda^* \y^*) &=& \sup_{\bq \in \Q} \{ \left< \bq, \Lambda^* \y^* \right> - F(\bq) \} =
                        \sup_{\bq \in \Q} \{ \left< \bq, \Lambda^* \y^* + \h \right> - \half \IntO  \B \bq \cdot \bq \dx \} \\
                        && (\mbox{use the constrain } \Lambda^* \by^*_2  + h_2 = -(\Lambda^* \by^*_1 + h_1))  \\
         &=& \sup_{(\q_1, \q_2) \in Q} \{ \left< \q_1 - \q_2, \Lambda^* \by^*_1 + h_1 \right> - \half \IntO \kappa(\q_1 - \q_2)^2 \dx \} \\
         && (\mbox{supremum is obtained for } \q_1 - \q_2 = (\Lambda^* \by^*_1 + h_1)/\kappa) \\
         &=&    \frac{1}{2 \kappa} \IntO (\Lambda^* \by^*_1 + h_1)^2 \dx = \frac{1}{4 \kappa}
\IntO \left[ (\Lambda^* \by^*_1 + h_1)^2 + (\Lambda^* \by^*_2 + h_2)^2 \right] \dx \\ &=&  \frac{1}{4 \kappa} \IntO (\Lambda^* \y^* + \h)^2 \dx.
\ee
\end{proof}
\end{lemma}
\begin{remark}
We note that \eqref{subspace} is a weaker restriction than the sum of two equilibrium relations
$\Lambda^* \by^*_1 + h_1 =0$ and $\Lambda^* \by^*_2 + h_2 =0$, which one would await from the general theory.
In other words, our analysis shows that the strict equilibrium of the dual variables in the componentwise sense is
not required in the couple system.
\end{remark}
After the substitution of \eqref{G_dual} and \eqref{F_dual} in the definition \eqref{compound_functionals}, we obtain
explicit expressions for the compound functionals.
\begin{lemma}[compound functionals] It holds
\ben
D_G(\Lambda \bq, \y^*)
&=& \half \IntO \A (\Lambda \bq - \A^{-1} \y^*): (\Lambda \bq - \A^{-1} \y^*) \dx \label{G_compound}, \\
D_F(\bq,\Lambda^* \y^*) &=& \left\{ \begin{array}{ll} \half \IntO  \B \bq \cdot \bq \dx
+ \frac{1}{4 \kappa} \IntO (\Lambda^* \y^* + \h)^2 \dx \cr
   \quad \quad \mbox{if } \Lambda^* \by^*_1 + h_1  + \Lambda^* \by^*_2 + h_2 = 0, \cr
+ \infty \quad \mbox{otherwise}.
\end{array}\right. \label{F_compound}
\een
\label{lemma_compounds}

\end{lemma}
According to \eqref{estimate_majorant}, the sharpest bound of $a(\bp-\bq, \bp - \bq)$ is provided by the estimate
\ben
\half a(\bp-\bq, \bp - \bq) \leq \inf_{\y^* \in \Y^*} M(\bq,\y^*).
\label{majorant_infimum}
\een
Since $M(\bq,\y^*)=+\infty$ if $\y^*$ does not satisfy \eqref{subspace}, we must restrict ourselves to arguments $\y^* \in \Yh$, where
\ben
\Yh:= \{(y^*_1, y^*_2) \in \Y^* : \Lambda^* \by^*_1 + h_1  + \Lambda^* \by^*_2 + h_2 = 0 \ae \}.
\label{constrain}
\een
To construct an element of $\Yh$ requires
an exact equilibration procedure, which have been studied for a Poisson problem in \cite{BrSc}.
Below, we show a way to avoid the constrain \eqref{constrain} by a special penalty term added to the
functional majorant. We define
\ben
\Ydiv:= \{(\by^*_1, \by^*_2) \in \Y^* : \Lambda^* \by^*_1 + \Lambda^* \by^*_2 \in L^2(\Omega)  \}
\een
and note that $ \Yh \subset \Ydiv $ (since $h_1, h_2 \in L^2(\Omega)$).
Further we decompose
$$\y^* =  \hat \y^* + (\y^* - \hat \y^*) $$
with $\hat \y^* \in \Ydiv $ and we extend the dual functionals $D_G$ and $D_F$ by the new variable $\hat \y^*$. We rewrite \eqref{G_compound} as
\be
&&D_G(\Lambda \bq, \y^*) = \half \IntO \A (\Lambda \bq - \A^{-1} \hat \y^*): (\Lambda \bq - \A^{-1} \hat \y^*) \dx +\\
&&+ \IntO (\Lambda \bq - \A^{-1} \hat \y^*): (\y^* - \hat \y^*) \dx
+ \half \IntO \A^{-1} (\y^* - \hat \y^*): (\y^* - \hat \y^*) \dx
\ee
and use the inequality $2 \M_1:\M_2 \leq \beta_1 M_1:M_1 + \frac{1}{\beta_1} \M_2:\M_2$ valid for all matrices $\M_1, \M_2$ and for all $\beta_1 > 0$ to bound
the middle term as
\begin{multline}
(\Lambda \bq - \A^{-1} \hat \y^*): (\y^* - \hat \y^*) = \A^{1/2} (\Lambda \bq - \A^{-1} \hat \y^*): \A^{-1/2} (\y^* - \hat \y^*) \\
\leq \frac{\beta_1}{2} \A (\Lambda \bq - \A^{-1} \hat \y^*):(\Lambda \bq - \A^{-1} \hat \y^*) + \frac{1}{2 \beta_1} \A^{-1} (\y^* - \hat \y^*): (\y^* - \hat \y^*).
\end{multline}
Obviously, the middle terms adds to the left and the right terms in $D_G(\Lambda \bq, \y^*)$ above and the modified compound functional reads
\begin{equation}
\begin{split}
D_G(\Lambda \bq, \y^*, \hat \y^*) :=& \frac{1+\beta_1}{2} \IntO \A (\Lambda \bq - \A^{-1} \hat \y^*): (\Lambda \bq - \A^{-1} \hat \y^*) \dx \\
&+ (\half + \frac{1}{2 \beta_1}) \IntO \A^{-1} (\y^* - \hat \y^*): (\y^* - \hat \y^*) \dx. \label{G_compount_modified}
\end{split}
\end{equation}
It also contains a scalar factor $\beta_1 > 0$ that value can be chosen arbitrarily.
Similar technique is used to modify the compound functional $D_F(\bq,\Lambda^* \y^*)$.
For the second integral in \eqref{F_compound}, we have
\be
\IntO (\Lambda^* \y^* + \h)^2 \dx \leq (1 + \beta_2) \IntO (\Lambda^* \hat \y^* + \h)^2 \dx + (1 + \frac{1}{\beta_2}) \IntO (\Lambda^* (\y^* - \hat \y^* ))^2 \dx,
\ee
where $\beta_2 > 0$. Therefore, a modified dual functional reads
\begin{equation}
\begin{split}
D_F(\bq,\Lambda^* \y^*, \Lambda^* \hat \y^*) :=&
\half \IntO  \B \bq \cdot \bq \dx + \frac{1}{4 \kappa} (1 + \beta_2) \IntO (\Lambda^* \hat \y^* + \h)^2 \dx \\
&+\frac{1}{4 \kappa} (1 + \frac{1}{\beta_2}) \IntO (\Lambda^* (\y^* - \hat \y^* ))^2 \dx \label{F_compount_modified}.
\end{split}
\end{equation}
By adding \eqref{G_compount_modified} and \eqref{F_compount_modified}, we extend the functional majorant \eqref{majorant} to
\ben M(\bq,\y^*, \hat \y^*) := D_F(\bq,\Lambda^* \y^*, \Lambda^* \hat \y^*) + D_G(\Lambda \bq, \y^*, \hat \y^*), \een
in which arbitrary variables satisfy the constrain
$$(\y^*, \hat \y^*) \in \Yh \times \Ydiv.$$
Clearly, the original and extended majorants satisfy the inequality
\ben
\half a(\bp-\bq, \bp - \bq) \leq M(\bq,\y^*) \leq M(\bq,\y^*, \hat \y^*) \label{estimate_general}
\een
for all $\hat \y^* \in \Ydiv, \beta_1 > 0, \beta_2 >0$. This estimate is sharp in the sense that there are no irremovable gaps in the inequalities.
 Indeed, if we set $\y^*= \hat \y^* = \Lambda \p$ and tend $\beta_1$ and $\beta_2$ to zero,
then $M(\bq,\y^*, \hat \y^*)$ tends to $M(\bq,\y^*)$ (and even to the exact error $\half a(\bp-\bq, \bp - \bq)$,
cf. \eqref{majorant_infimum}).

\subsection{An upper estimate of $M(\bq,\y^*, \hat \y^*)$}
Let us denote $\y^*=(\by^*_1, \by^*_2)$ and $\hat \y^*=(\hat \by^*_1, \hat \by^*_2)$ and consider a particular subspace
\ben
(\y^*, \hat \y^*) \in \{ \Yh \times \Ydiv:  \Lambda^* \by^*_1 + h_1  = 0, \by^*_2=\hat \by^*_2 \enspace \ae \} \label{subspace_choice1}.
\een
In this subspace, it holds (cf. \eqref{lambda_dual})
\be
\IntO (\Lambda^* (\y^* - \hat \y^* ))^2 \dx = \IntO (\diver (\hat \by^*_1 - \by^*_1 ))^2 \dx =
\IntO (\diver \hat \by^*_1 - h_1)^2 \dx.
\ee
Therefore, $D_F(\bq,\Lambda^* \y^*, \Lambda^* \hat \y^*)$ defined in \eqref{F_compount_modified} simplifies as $\y^*$-independent
\ben
D_F(\bq, \Lambda^* \hat \y^*) &:=&
\half \IntO  \B \bq \cdot \bq \dx + \frac{1}{4 \kappa} (1 + \beta_2) \IntO (\Lambda^* \hat \y^* + \h)^2 \dx \label{DF_final} \\
&&+\frac{1}{4 \kappa} (1 + \frac{1}{\beta_2}) \IntO (\diver \hat \by^*_1 - h_1)^2 \dx \nonumber
\een
and only $\y^*$-dependent functional in $D_G(\Lambda \bq, \y^*, \hat \y^*)$ defined in \eqref{G_compount_modified} writes
\ben
\IntO \A^{-1} (\y^* - \hat \y^*): (\y^* - \hat \y^*) \dx = \IntO k_1^{-1} (\by^*_1 - \hat \by^*_1)\cdot (\by^*_1 - \hat \by^*_1) \dx.
\label{intermediate}
\een
\begin{lemma} \label{lemma1}
Let us  define a space
$$
Y_{h_1}:= \{\by^*_1 \in L^2(\Omega)^d : \Lambda^* \by^*_1 + h_1  = 0 \ae \}.
$$
Then, for all $\hat \by^*_1 \in \Hdiv$, it holds
$$\inf_{\by^*_1 \in Y_{h_1}} \IntO \Norm{\by^*_1 - \hat \by^*_1}^2 \dx \leq C^2 \Norm{\diver \hat \by^*_1 + h_1 }^2$$
where $C > 0$ satisfies Friedrichs' inequality
$\Norm{w}_{L^2(\Omega)} \leq C \Norms{\nabla w}{L^2(\Omega)}$ valid for all  $w \in H^1_0(\Omega)$.
\begin{proof} It follow from Theorem 6.1 from \cite{ReVa2} by the modification related to the fact the we
consider vector arguments.
\end{proof}
\end{lemma}
Application of Lemma \ref{lemma1} to \eqref{intermediate} and the back substitution to \eqref{G_compount_modified} defines a $\y^*$-independent dual functional
\begin{multline}
D_G(\Lambda \bq, \hat \y^*) := \frac{1 + \beta_1}{2} \IntO \A (\Lambda \bq - \A^{-1} \hat \y^*): (\Lambda \bq - \A^{-1} \hat \y^*) \dx \\
+ k_1^{-1} (\half + \frac{1}{2 \beta_1}) C^2 \Norm{\diver \hat \by^*_1 + h_1 }^2. \label{DG_final}
\end{multline}
which provides an upper estimate of the quantity
$$\inf_{\y^* \in \Yh} D_G(\Lambda \bq, \y^*, \hat \y^*).$$
Therefore, the sum of \eqref{DF_final} and \eqref{DG_final} defines a $\y^*$-independent functional
\ben M_{\beta_1, \beta_2}(\bq, \hat \y^*) := D_F(\bq, \Lambda^* \hat \y^*) + D_G(\Lambda \bq, \hat \y^*) \een
that serves as an upper bound of $M(\bq,\y^*, \hat \y^*)$ and provides a computable estimate
\ben
\half a(\bp-\bq, \bp - \bq) \leq M_{\beta_1, \beta_2}(\bq, \hat \y^*) \quad \mbox{for all } \hat \y^* \in \Ydiv.
\label{estimate_majorant}
\een
\begin{remark}[symmetric form of $D_G$]
If we replace the subspace \eqref{subspace_choice1} by
\ben
(\y^*, \hat \y^*) \in \{ \Yh \times \Ydiv:  \Lambda^* \by^*_2 + h_2  = 0, \by^*_1=\hat \by^*_1 \enspace \ae \},
\een
then, instead of \eqref{DG_final}, we obtain
\begin{multline}
D_G(\Lambda \bq, \hat \y^*) := \frac{1 + \beta_1}{2} \IntO \A (\Lambda \bq - \A^{-1} \hat \y^*): (\Lambda \bq - \A^{-1} \hat \y^*) \dx \\
+ k_2^{-1} (\half + \frac{1}{2 \beta_1}) C^2 \Norm{\diver \hat \by^*_2 + h_2 }^2 \label{DG_final2}.
\end{multline}

\end{remark}

\section{A posteriori error estimate for approximations of the coupled system (1.1)}

Assume that the fluid pressures $\p_1$ and $\p_2$ are resolved exactly and substituted to the elasticity equation (cf. ~\eqref{system})
\ben
 -\nabla \cdot (\bbL \eps(u)) = \f(x,t) + \alpha_1 \nabla \p_1 + \alpha_2 \nabla \p_2. \nonumber
\een
Let $\v$ be an approximation of $\u$ (this problem is considered in the same domain $\Omega$ as the problem
 \eqref{diffusion1}-\eqref{diffusion2}). We define the Dirichlet boundary condition by a function
$\u_0 \in H^1(\Omega; \Rd)$ and assume
$$\v\in \u_0 + H^1_0(\Omega; \Rd).$$
\begin{lemma}
For every function
$\ttau \in Q :=\{\sigma \in L^2(\Omega; \Mdsym): \diver \sigma \in L^2(\Omega; \Rd)\}$
it holds
\begin{multline}
\Norm{\eps(\u-\v)}_{\bbL; \Omega} \leq
 \Norm{\eps(\v) - \bbL^{-1} \ttau}_{\bbL; \Omega} + C \Norm{\diver \ttau + \f - \alpha_1 \nabla \p_1 - \alpha_2 \nabla \p_2}_{\Omega},
\label{estimate_elastic}
\end{multline}
where the constant $C > 0 $ satisfies an inequality
\ben
\label{Korns}
\Norm{\w}_{\Omega} \leq C \Norm{\eps(\w)}_{\bbL; \Omega} \quad \mbox{for all } \w \in H^1_0(\Omega; \Rd).
\een
and the norm $\Norm{\cdot}$ is defined as
$$
\Norm{\eps}_{\bbL; \Omega}^2 :=\IntO \bbL \eps : \eps \dx.
$$
\begin{proof} Estimates in chapter 6.5 in \cite{Re3} which are applied to the linear elasticity problem
with the right-hand side $\f - \alpha_1 \nabla \p_1 - \alpha_2 \nabla \p_2$.
The existence of constant $C$ follows from Korn's and Friedrichs' inequalities.
\end{proof}
\end{lemma}
\begin{remark}
The estimate ~\eqref{estimate_elastic} is sharp with respect to parameter $\tau$. Indeed, the choice $\tau= \bbL \eps(u)$ satisfies the equilibrium condition
\ben
\diver \tau + \f = \alpha_1 \nabla \p_1 + \alpha_2 \nabla \p_2
\label{equilibrium_elastic}
\een
and reduces therefore ~\eqref{estimate_elastic} to the equality.
\end{remark}
Let $\q_1$ and $\q_2$ be approximation of exact pressure fields $p_1$ and $p_2$ respectively.
By triangle inequalities, we obtain
\begin{multline}
\Norm{\diver \tau + \f - \alpha_1 \nabla \p_1 - \alpha_2 \nabla \p_2}_{\Omega} \leq \Norm{\diver \tau + \f - \alpha_1 \nabla \q_1 - \alpha_2 \nabla \q_2}_{\Omega} \\
+ \Norm{\nabla (\p_1 - \q_1)}_{\Omega}  + \Norm{\nabla (\p_2 - \q_2)}_{\Omega}. \label{4.4}
\end{multline}
Use \eqref{4.4} and square both parts of \eqref{estimate_elastic} to obtain
\ben
\Norm{\eps(\u-\v)}_{\bbL; \Omega}^2 &\leq &
 ( \Norm{\eps(\v) - \bbL^{-1} \ttau}_{\bbL; \Omega} \nonumber \\
&&+ C \Norm{\diver \ttau + \f - \alpha_1 \nabla \q_1 - \alpha_2 \nabla \q_2}_{\Omega}  \label{4.5}\\
 && + C \Norm{\nabla (\p_1 - \q_1)}_{\Omega} + C \Norm{\nabla (\p_2 - \q_2)}_{\Omega} )^2 \nonumber .
\een
By the algebraic inequality
$$(a + b + c)^2 \leq (1 + \beta_4 + \beta_5)\enspace a^2 + (1+ \frac{1}{\beta_4} + \beta_6) \enspace b^2 + (1+ \frac{1}{\beta_5} + \frac{1}{\beta_6}) \enspace c^2$$
valid for all scalars $a, b, c$ and for all $\beta_4, \beta_5, \beta_6 > 0$, inequality \eqref{4.5}
and the following inequality ($\beta_3$ is an arbitrary positive constant)
\begin{multline}
\left(\Norm{\nabla (\p_1 - \q_1)}_{\Omega}  + \Norm{\nabla (\p_2 - \q_2)}_{\Omega}\right)^2 \\
\leq (1 + \beta_3) \Norm{\nabla (\p_1 - \q_1)}_{\Omega}^2  + (1 + \frac{1}{\beta_3}) \Norm{\nabla (\p_2 - \q_2)}_{\Omega}^2 \\
\leq \max \{ \frac{1 + \beta_3}{k_1}, \frac{1 + \beta_3}{k_2 \beta_3}  \}
\enspace a(\bp-\bq, \bp - \bq) \\
\leq 2 \max \{ \frac{1 + \beta_3}{k_1}, \frac{1 + \beta_3}{k_2
\beta_3}  \} \enspace M_{\beta_1, \beta_2}(\bq, \hat \y^*).
\end{multline}
Now  we obtain the final estimate in terms of the coupled
error norm
\begin{multline}
a(\bp-\bq, \bp - \bq)+\Norm{\eps(\u-\v)}_{\bbL; \Omega}^2
 \leq  (1+\beta_4 + \beta_5) \Norm{\eps(\v) - \bbL^{-1} \ttau}_{\bbL; \Omega}^2 \\
 + \left(1+ \frac{1}{\beta_4} + \beta_6\right) C^2
 \Norm{\diver \ttau + \f - \alpha_1 \nabla \q_1 - \alpha_2 \nabla \q_2}_{\Omega}^2
 +2 \widehat C
  \enspace M_{\beta_1, \beta_2}(\bq, \hat \y^*),
\label{estimate_elastic_next}
\end{multline}
where
\be
\widehat C=1+ C^2\left(1+ \frac{1}{\beta_5} + \frac{1}{\beta_6}\right)  \max
\left\{ \frac{1 + \beta_3}{k_1}, \frac{1 + \beta_3}{k_2 \beta_3}  \right\}.
\ee
This estimate
holds for all $\ttau \in Q, \hat \y^* \in \Ydiv$ and all $\beta_1, \dots, \beta_6 >0$.

%
%


\end{document}